%January 23, 2007
\input amstex
\documentstyle{amsppt}
\NoBlackBoxes
\TagsOnRight
\NoRunningHeads
\hcorrection{2cm}

\loadbold\loadeusm\loadeufm
\def\qed{\hphantom{a}\hfill $\square$}
\def\nin{\notin}
\def\Ra{\Rightarrow}
\def\lb{\linebreak}

\def\mb{\mathbreak}\def\nmb{\nomathbreak}

\def\supp{\operatorname{supp}}

\def\Aut{\operatorname{Aut}}
\def\Inn{\operatorname{Inn}}

\def\<{\langle}\def\>{\rangle}

\def\wlim{\operatorname{w-lim}}
\def\w*lim{\operatorname{w*-lim}}

\def\vep{\varepsilon}\def\vp{\varphi}

\def\sub{\subseteq}\def\bus{\supseteq}

\def\ov{\overline}

\def\G{\varGamma}\def\g{\gamma}
\def\om{\omega}
\def\Om{\varOmega}
\def\al{\alpha}

\def\lam{\lambda}\def\del{\delta}
\def\iy{\infty}
\def\si{\sigma}
\def\k{\kappa}
\def\ti{\times}

\def\inl{\int\limits}

\topmatter
\title{The Choquet-Deny theorem and distal properties of totally 
disconnected locally compact groups of polynomial growth} 
\endtitle
\author{\smc Wojciech Jaworski* and C. Robinson Edward Raja**}
\vskip1pt {\rm School of Mathematics and 
Statistics \\ Carleton University\\
Ottawa, Ontario, Canada K1S 5B6}
\vskip2pt
{\rm e-mails:\ {\it wjaworsk\@math.carleton.ca\/ {\rm and} 
creraja\@isibang.ac.in\/}}
\endauthor 
%\leftheadtext{}
%\rightheadtext{}
%
\abstract We obtain sufficient and necessary conditions for 
the Choquet-Deny theorem to hold in the class of compactly generated totally 
disconnected locally compact groups of polynomial growth, and in a larger 
class of totally disconnected generalized $\ov{FC}$-groups. The following 
conditions turn out to be equivalent when $G$ is a metrizable compactly 
generated totally disconnected locally compact group of polynomial growth:
(i) the Choquet-Deny theorem holds for $G\/$; (ii) the group of inner 
automorphisms of $G$ acts distally on $G\/$; (iii) every inner automorphism 
of $G\/$ is distal; (iv) the contraction subgroup of every inner automorphism 
of $G\/$ is trivial; (v) $G$ is a SIN group. We also show that for every 
probability measure $\mu$ on a totally disconnected compactly generated 
locally compact second countable group of polynomial 
growth, the Poisson boundary is a homogeneous space of $G$, and that it is 
a compact homogeneous space when the support of $\mu$ generates $G$. 
\medpagebreak

\noindent {\it 2000 Mathematics Subject Classification}\/: 60B15, 43A05, 
60J50, 22D05, 22D45
\endabstract  
\endtopmatter
\document
\footnote""{\hphantom{*}*\ Supported by an NSERC Grant.}
\footnote""{**\ Permament address: Indian Statistical Institute, Statistics 
and Mathematics Unit, 8th Mile Mysore Road, Bangalore 560 059, India.}
\head{\bf 1. Introduction}\endhead

Let $\mu$ be a regular Borel probability measure on a locally compact group 
$G$. A bounded 
Borel function $h\:G\to\Bbb C$ is called $\mu$-{\it harmonic\/} if it 
satisfies 
$$
h(g)=\inl_G h(gg')\,\mu(dg'), \quad g\in G.\tag{1.1}
$$ 
The classical Choquet-Deny theorem asserts that when $G$ is abelian then 
every bounded continuous $\mu$-harmonic function is constant on the 
(left) cosets of the smallest closed subgroup, $G_\mu$, containing the 
support of $\mu$. 

The Choquet-Deny theorem remains true for many nonabelian locally compact 
groups, e.g., 2-step nilpotent groups [10], nilpotent 
[SIN] groups [11], and compact groups. But it does not hold for all groups. 
If the theorem holds for a probability measure $\mu$ then $G_\mu$ must 
necessarily be an amenable subgroup [7,27]. It follows that groups for 
which the theorem is valid are necessarily amenable. However, the theorem is 
not true for every amenable group [20]. 

The stronger condition, that $G$ have polynomial growth, is sufficient for 
the theorem to hold when $G$ is a finitely generated (discrete) group 
[20,18]. When $G$ is finitely generated and solvable then the theorem holds 
if and only if $G$ has polynomial growth [18]. In general, the theorem 
fails for discrete groups of polynomial growth that are not finitely 
generated, in particular, it is not true for locally finite groups [20]. 
It appears that the largest class of discrete groups known today for which 
the Choquet-Deny theorem is true is the class of FC-hypercentral groups [12]. 
This class is a proper subclass of the class of discrete groups of 
polynomial growth, while finitely generated FC-hypercentral groups are 
precisely the finitely generated groups of polynomial growth. We do not know 
of any discrete groups for which the Choquet-Deny theorem is true and which 
are not FC-hypercentral. 

A probability measure $\mu$ on a locally compact group $G$ is called {\it 
spread out\/} if for some $n$ the convolution power $\mu^n$ is nonsingular. 
With the restriction that $\mu$ be spread out the Choquet-Deny theorem holds 
for all locally compact nilpotent groups [2,16] and for compactly generated 
locally compact groups of polynomial growth [16]. When $G$ is almost 
connected, then $G$ has polynomial growth if and only if the Choquet-Deny 
theorem holds for every spread out measure [16]. The same is true when $G$ is 
a Zariski-connected $p$-adic algebraic group [24, Theorem 4.2]. 
While it remains an open question whether the spread out assumption can be 
disposed of when $G$ is nilpotent\,
\footnote{\ The recently published proof [25] is incomplete (the proof of 
Lemma 2.5 has a gap).}
, it is known that the Choquet-Deny theorem is not true for arbitrary 
probability measures on compactly generated locally compact groups of 
polynomial growth [13, Remark 3.15]. 

The main goal of the present article is to obtain necessary and sufficient 
conditions for the validity of the Choquet-Deny theorem in the class of 
compactly generated totally disconnected locally compact groups of polynomial 
growth, and in a larger class of totally disconnected `generalized 
$\ov{FC}$-groups' [3,21]. It turns out that the key to finding such 
conditions is a study of distal properties of totally disconnected groups. 
This motivates our investigations in the next section, which can also be of 
quite independent interest. The Choquet-Deny theorem for generalized 
$\ov{FC}$-groups is discussed in Section 3. In Section 4 we remark on the 
structure of boundaries of random walks on compactly generated totally 
disconnected groups of polynomial growth and certain related groups. 
\head{\bf 2. Distal properties of totally disconnected locally compact 
groups}\endhead

Let $G$ be a Hausdorff topological group and $\G$ a subgroup of $\Aut(G)$, 
the group of topological automorphisms of $G$. We will say that $\G$ is 
{\it distal\/} (or acts distally on $G$) if for any $x\in G-\{e\}$, the 
identity element $e\/$ is not in the closure of the orbit $\G x=\{\g 
(x)\,;\,\g\in \G\}$. A single automorphism $\g\in\Aut(G)$ will be called 
distal if the subgroup $\<\g\>$ it generates acts distally on $G$. An element 
$g$ of $G$ will be called distal if the corresponding inner automorphism 
$\g(\cdot)=g\cdot g^{-1}$ is distal. We will say that $G$ is distal if the 
group $\Inn(G)$ of inner automorphisms of $G$ acts distally 
on~$G$. 

Trivially, if $G$ is distal then every $g\in G$ is distal. While the converse 
is not true in general, Rosenblatt [26] proved that when $G$ is an almost 
connected locally compact group then $G$ is distal if and only if every $g\in 
G$ is distal; moreover $G$ is distal if and only if it has polynomial growth. 
According to [23] this remains true also for certain classes of $p$-adic Lie 
groups. However, there are many locally compact 
groups of polynomial growth that are not distal. For example, the semidirect 
product $K\ti\negmedspace_\tau\medspace\Bbb Z$ where $K$ is a nontrivial 
compact metric group and $\tau$ is an ergodic automorphism of $K$, will never 
be distal. 

Given $\g\in\Aut(G)$ the {\it contraction subgroup\/} of $\g$ is the subgroup 
$C(\g)=\{x\in G\,;\,\mb\lim_{n\to\iy}\g^n(x)= e\}$. When $\g$ is the inner 
automorphism $\g(\cdot)=g\cdot g^{-1}$, we will write $C(g)$ for $C(\g)$. 
Obviously, if $\tau\in\Aut(G)$ is distal then $C(\tau)=C(\tau^{-1})=\{e\}$. 
When $G$ is a Lie group, the three conditions: {\sl $\G$ is distal\/}; 
{\sl every $\g\in\G$ is distal\/}; and, {\sl $C(\g)=\{e\}$ for every 
$\g\in\G\/$}, are equivalent for every subgroup $\G$ of $\Aut(G)$  [1]. 

Recall that a subgroup $\G$ of $\Aut(G)$ is equicontinuous (at $e$) if and 
only if $G$ admits a neighbourhood base at $e$ consisting of neighbourhoods 
that are invariant under $\G$. When $G$ is locally compact and totally 
disconnected then $\G\/$ is equicontinuous if and only if compact open 
subgroups invariant under $\G\/$ form a neighbourhood base at $e$. 
Equicontinuous automorphism groups are obviously distal. A SIN group is a 
topological group $G$ for which $\Inn(G)$ is equicontinuous. SIN groups 
are distal but, in general, distal groups are not SIN groups (e.g., a 
nilpotent group need not be SIN but every nilpotent group is distal [26]). 

Our goal in this section is to prove that for a class of compactly generated 
totally disconnected locally compact groups, including groups of polynomial 
growth, the four conditions: 
{\sl $G$ is distal\/}; {\sl every $g\in G$ is distal\/}; {\sl $C(g)=\{e\}$ 
for every $g\in G$\/}; and, {\sl $G$ is SIN\/}, are equivalent. In the 
following section we will show that for this class of groups the four 
conditions and the condition that $G$ have polynomial growth, are equivalent 
to the condition that the Choquet-Deny theorem hold for $G$. 

Some of the recent results of Baumgartner and Willis on contractions 
subgroups [4], based on Willis' theory of tidy subgroups [29], play a 
key role in our argument. These results are proven for metrizable groups, 
hence, in many of our results we need to assume metrizability.  

\proclaim{Proposition 2.1} If $G$ is a totally disconnected metrizable 
locally compact group then for every $\tau\in\Aut(G)$ the following 
conditions are equivalent\rom: 
\roster
\item"{(i)}"$\tau$ is distal,
\item"{(ii)}"$C(\tau)=C(\tau^{-1})=\{e\}$,
\item"{(iii)}"for every compact open subgroup $U$ there exists $k=0,1,\dots$ 
such that $\tau\bigl(\bigcap_{i=0}^k\tau^i(U)\bigr) 
=\bigcap_{i=0}^k\tau^i(U)$,
\vskip2pt
\item"{(iv)}"$\<\tau\>$ is equicontinuous. 
\endroster
\endproclaim
\demo{Proof}The only nonobvious implication in the chain 
(i)$\Ra$(ii)$\Ra$(iii)$\Ra$(iv)$\Ra$(i) is (ii)$\Ra$(iii). Let $U$ be a 
compact open subgroup. Since $C(\tau)$ is closed, by [4, Theorem~3.32] 
there exists $k$ such that $V= \bigcap_{i=0}^k\tau^i(U)$ is tidy for $\tau$. 
But as $C(\tau)=C(\tau^{-1})=\{e\}$, [4, Proposition 3.24] implies that 
$s(\tau)=s(\tau^{-1})=1$  where $s\:\Aut(G)\to\Bbb N$ is the scale function. 
Since $s(\tau)=[\tau (V):V\cap\tau (V)]$ and 
$s(\tau^{-1})=[\tau^{-1}(V):V\cap\tau^{-1}(V)]$, so $\tau(V)=V$. \qed\enddemo

\proclaim{Lemma 2.2} Let $\G$ be a subgroup of $\Aut(G)$ where $G$ is a 
totally disconnected metrizable locally compact group. If 
$\tau_1,\tau_2,\dots,\tau_n\in\Aut(G)$ are distal and 
for every $j=1,2,\dots,n$, $\,[\tau_j,\<\G\cup\{\tau_1,\dots,\tau_{j-1}\}\>] 
\sub \<\G\cup\{\tau_1,\dots,\tau_{j-1}\}\>$, then for every compact open 
subgroup $U$ invariant under $\G$ there exists a compact open subgroup 
$V\sub U$ invariant under $\<\G\cup\{\tau_1,\dots,\tau_n\}\>$.
\endproclaim

\demo{Proof}It is clear that the lemma follows by induction once 
it is verified for $n=1$. So we suppose that $n=1$. 

By Proposition 2.1 there exists $k$ such that $V=\bigcap_{i=0}^k\tau^i_1(U)$ 
satisfies $\tau_1(V)=V$. It is enough to show that $\g(V)=V$ for every $\g\in 
\G$. But our assumption implies that $[\tau^i_1,\G]\sub\G$ for every 
$i=0,1,\dots$\,. Hence, given $\g\in\G$ we obtain 
$\g(V)=\bigcap_{i=0}^k(\g\tau^i_1)(U)=
\bigcap_{i=0}^k(\tau^i_1\g[\g,\tau^i_1])(U) =\bigcap_{i=0}^k\tau^i_1(U)=V$. 
\qed\enddemo
\proclaim{Lemma 2.3}Let $\G$ be a subgroup of $\Aut(G)$ where $G$ is a 
totally disconnected metrizable locally compact group. Suppose that every 
$\g\in\G$ is distal and that $\G$ has a normal equicontinuous subgroup $\G_1$ 
with the quotient $\G/\G_1$ containing a polycyclic subgroup of finite index. 
Then $\G$ is equicontinuous. 
\endproclaim
\demo{Proof}Let $\Om$ be a neighbourhood of $e$. Denote by $P$ the polycyclic 
subgroup of finite index in $\G/\G_1$ and let $P_0=P$, $P_1=[P,P]$, 
$P_2=[P_1,P_1],\dots, P_m=\{\G_1\}$ be the derived series for $P$. Write 
$\pi$ for the canonical homomorphism $\pi\:\G\to\G/\G_1$ and put $\hat 
P_j=\pi^{-1}(P_j)$ for $j=0,1,...,m$. 

Suppose that for some $j=1,2,...,m$, $V\sub \Om$ is a compact open subgroup 
invariant under $\hat P_j$. We will show that there is 
then a compact open subgroup $W\sub V$ invariant under $\hat P_{j-1}$. Now, 
since $P$ is polycyclic, $P_{j-1}$ is generated 
by a finite set $\{p_1,\dots, p_n\}$. For every $i=1,2,\dots,n$ find 
$\tau_i\in\hat P_{j-1}$ with $p_i=\pi(\tau_i)$. Applying Lemma 2.2 to $\hat 
P_j$ and $\tau_1,\dots,\tau_n$ we conclude that there is a compact open 
subgroup $W\sub V$ invariant under $\<\hat 
P_j\cup\{\tau_1,\dots,\tau_n\}\>=\hat P_{j-1}$. 

Our assumption is that there is a compact open subgroup $V\sub \Om$ such that 
$\g(V)=V$ for every $\g\in \hat P_m=\G_1$. With the aid of the preceding 
paragraph it then follows that there is a compact open subgroup $W\sub\Om$ 
invariant under $\hat P_0$. Next, 
since $P=P_0$ has finite index in $\G/\G_1$, $\hat P_0$ has finite index in 
$\G$. Hence, the intersection $U=\bigcap_{\g\in\G}\g(W)$ is a compact open 
subgroup invariant under $\G$ and contained in $\Om$. 
\qed\enddemo

\proclaim{Corollary 2.4} Let $G$ be a totally disconnected metrizable locally 
compact group. If a subgroup $\G$ of $\Aut(G)$ contains a polycyclic subgroup 
of finite index then the following conditions are equivalent\rom:
\roster
\item"{(i)}"$\G$ is distal,
\item"{(ii)}"every $\g\in\G$ is distal,
\item"{(iii)}"$\G$ is equicontinuous. 
\endroster
\endproclaim 

As the following examples show, `polycyclic' in Corollary 2.4 cannot be 
replaced by `solvable'. In fact, the three conditions can be different for 
countable abelian groups of automorphisms. We do not know if `polycyclic' can 
be replaced by `finitely generated solvable'. 

\example{Example 2.5} Let $\vp\:\Bbb R\to\Bbb T$ denote the function 
$\vp(t)=e^{2\pi it}$ and let $H$ be any infinite subgroup of $\vp(\Bbb Q)$. 
Note that every $h\in H$ has finite order. Let $G$ be the totally 
disconnected compact abelian group $G=\Bbb Z_2^H$. $H$ acts on $G$ by left 
translations: $(hf)(x)=f(h^{-1}x)$ ($h\in H$, $f\in G$, $x\in H$). Let $\G$ 
be the resulting subgroup of $\Aut(G)$. Then every 
element of $\G$ is distal because it has finite order. However, $\G$ is not 
distal. Indeed, let $f\in G$ be the function $f(x)=\del_{1x}$ and let $U$ be 
any neighbourhood of $e$ in $G$. Then for some finite subset $F\sub H$, $U$ 
contains the set $\{g\in G\,;\,g(x)=0 \text{ for every }x\in F\}$. Hence, if 
$h\in H-F$ then $(hf)(x)=0$ for every $x\in F$, i.e., $hf\in U$. 

Thus for a countable abelian group of automorphisms (ii) does not imply (i) 
(nor (iii)). 
\endexample 

\example{Example 2.6} Let for $j\in\Bbb Z$, $G_j=\{x\in\Bbb Z_2^{\Bbb Z}\,;\, 
x_i=0\text{ for every }i\le j\}$ and let $G=\bigcup_{j\in\Bbb Z}G_j$. There 
is a locally compact totally disconnected group topology on $G$ in which the 
subgroups $G_j$, $j\in \Bbb Z$, form a neighbourhood base at $e$ (and are 
compact open). 
Given $j\in\Bbb N$ define $\tau_j\in\Aut(G)$ by $\tau_j(x)=y$ where $y_i=x_i$ 
for $i\in\Bbb Z-\{\pm j\}$, $y_j=x_{-j}$, and $y_{-j}=x_j$. The subgroup $\G$ 
of $\Aut(G)$ generated by $\tau_j$, $j\in\Bbb N$, is then abelian and 
distal: if $x\ne e$ and $j$ is the smallest integer with $x_j\ne 0$, then 
$\tau(x)\notin G_{|j|}$ for any $\tau\in\G$. However, $\G$ is not 
equicontinuous because if it were, there would exist $k\ge 0$ with 
$\tau(x)\in G_0$ for every $x\in G_k$ and $\tau\in\G$. However, if 
$x=(\del_{k+1\,n})_{n\in\Bbb Z}$ then $x\in G_k$, but 
$(\tau_{k+1}(x))_{-k-1}=x_{k+1}=1$, so that $\tau_{k+1}(x)\notin 
G_0$.

Thus for a countable abelian group of automorphisms (i) does not imply (iii).
\endexample 

\proclaim{Theorem 2.7}Suppose that a totally disconnected metrizable locally 
compact group $G$ admits an open normal SIN subgroup $N$ such that the 
quotient $G/N$ contains a polycyclic subgroup of finite index. Then the 
following conditions are equivalent:
\roster
\item"{(i)}"$G$ is distal,
\item"{(ii)}"every $g\in G$ is distal,
\item"{(iii)}"$G$ is SIN. 
\endroster
\endproclaim 

\demo{Proof}(ii)$\Ra$(iii): Let $\G=\Inn(G)$ 
and $\al:G\to \G$ be the canonical homomorphism. Put $\G_1=\al(N)$. 
$\G_1$ is a normal subgroup of $\G\/$ and $\G/\G_1$ is a homomorphic image of 
$G/N$. Therefore $\G/\G_1$ contains a polycyclic subgroup of finite index. 
Since $N$ is open and SIN, it follows that $\G_1$ is equicontinuous. Thus 
Lemma 2.3 applies. 
\qed\enddemo

Following [3] and [21] we call a locally compact group a {\it generalized 
$\ov{FC}$-group\/} if $G$ has a series $G=G_0\bus G_1\bus\dots\bus G_n=\{e\}$ 
of closed normal subgroups such that for every $i=0,1,\dots,n-1$, 
$G_i/G_{i+1}$ is a compactly generated group with precompact conjugacy 
classes. 
Every compactly generated locally compact group of polynomial growth 
is a generalized $\ov{FC}$ group [21, Theorem 2]. 
Every closed subgroup of a generalized $\ov{FC}$ group is compactly 
generated [21, Proposition 2]. A locally compact solvable group $G$ is a 
generalized $\ov{FC}$ group if and only if each closed subgroup of $G$ is 
compactly generated [10, Th\'eor\`eme III.1]. Using Propositions 1 and 7(ii) 
in [21] it is straighforward to give the following characterization of 
totally disconnected generalized $\ov{FC}$ groups: 
\proclaim{Proposition 2.8} A totally disconnected locally compact group $G$ 
is a generalized $\ov{FC}$ group if and only if it admits a compact open 
normal subgroup $N$ with the quotient $G/N$ containing a polycyclic subgroup 
of finite index. 
\endproclaim 

\proclaim{Theorem 2.9}Conditions {\rm (i), (ii)}, and {\rm (iii)} of 
Theorem 2.7 are equivalent when $G$ is a totally disconnected generalized 
$\ov{FC}$-group. 
\endproclaim 

\demo{Proof} When $G$ is metrizable, this is a special case of Theorem 2.7. 
We need to show that the implication (ii)$\Ra$(iii) is also true when $G$ is 
not metrizable. 

Note that $G$ is necessarily $\si$-compact (as it is compactly generated). 
Let $U$ be a neighbourhood of $e$ contained in the subgroup $N$ of 
Proposition 2.8. Find a neighbourhood $V$ 
of $e$ with $V^2\sub U$. By [8, Theorem 8.7] $V$ contains a compact normal 
subgroup $K$ such that $G/K$ is metrizable. Let $\pi\:G\to G/K$ denote the 
canonical homomorphism. Since $N$ is compact, we can use the theorem stating 
that a factor of a distal flow is distal [5, Corollary 6.10, p.\,52] to 
conclude that the restriction of every inner automorhism of $G/K$ to $N/K$ is 
distal. As $N/K$ is open, every $g\in G/K$ is then distal. Hence, by Theorem 
2.7, $G/K$ is SIN. Thus $\pi(V)$ contains a compact open normal subgroup $W$. 
Then $\hat W=\pi^{-1}(W)\sub VK\sub U$ and $\hat W$ is a compact open 
normal subgroup of $G$. 
\qed\enddemo
Since nilpotent groups are distal, Theorem 2.9 implies that a totally 
disconnected compactly generated locally compact nilpotent group is a SIN 
group, a result due to Hofmann, Liukkonen, and Mislove [9]. 

We note that for totally disconnected groups of polynomial growth which are 
not compactly generated, conditions (i),(ii),(iii) are different. In fact, 
the equivalence fails already for metabelian groups of polynomial growth. 
Examples of totally disconnected 2-step nilpotent groups which are not 
SIN groups can be found in [9] and [28]. An example of a metabelian group of 
polynomial growth which satisfies (ii) but not (i) (nor (iii)) is also 
readily available:

\example{Example 2.10}Let $H$ be as in Example 2.5 and let $W$ be the 
complete wreath product $W=\Bbb Z_2\wr H$. Evidently, $W$ is not distal but 
every $w\in W$ has finite order, so is distal. 
\endexample 

In the remainder of this section we prove that conditions (i),(ii),(iii) of 
Theorem 2.7 are equivalent for every metrizable compactly generated totally 
disconnected locally compact metabelian group. 

\proclaim{Lemma 2.11} If a locally compact group $G$ contains a normal 
finitely generated subgroup $N$ and a compact set $K$ such that 
$KN=G$, then $G$ is a SIN group.
\endproclaim 
\demo{Proof}Recall that the group of automorphisms of a finitely generated 
group is countable. Since the centralizer $C_G(N)$ of $N$ in $G$ is the 
kernel of the homomorphism which maps $g\in G$ to the restriction of the 
inner automorphism $g\cdot g^{-1}$ to $N$, it follows that $G/C_G(N)$ is 
countable. As $C_G(N)$ is a closed subgroup and $G$ is of the second 
category, we conclude that $C_G(N)$ is open. 

Let $U$ be a neighbourhood of $e$. Put $V=U\cap C_G(N)$ and let 
$V'$ be a neighbourhood of $e$ with $g^{-1}V'g\sub V$ for 
every $g\in K$. Then $W=\bigcap_{g\in G}gVg^{-1}=\bigcap_{g\in K}gVg^{-1}\bus 
V'$. Thus $W$ is a neighbourhood of $e$, invariant under $\Inn(G)$ and 
contained in $U$. 
\qed\enddemo

\proclaim{Proposition 2.12}If a compactly generated totally disconnected 
locally compact group $G$ contains a closed cocompact normal SIN 
subgroup, then $G$ is a SIN group.
\endproclaim 
\demo{Proof}Let $N$ denote the closed cocompact normal SIN subgroup and let 
a compact open subgroup $U$ of $G$ be given. A routine argument shows that 
$\Inn(G)$ acts equicontinuously on $N$. Hence, $U$ contains a compact 
subgroup $V$ of $N$ which is open in $N$ and normal in $G$. 

Let $\pi\:G\to G/V$ denote the canonical 
homomorphism. Since $N$ is cocompact, it is compactly generated [22]. 
Since $V$ is open in $N$, $\pi(N)$ is then finitely generated. It is also 
normal and there is a compact $K\sub G/V$ with $K\pi(N)=G/V$. 
Hence, by Lemma 2.11 $G/V$ is a (totally disconnected) SIN group. Thus 
$\pi(U)$ contains a compact open normal subgroup $W$. Then $\pi^{-1}(W)$ is a 
compact open normal subgroup contained in $UV=U$. 
\qed\enddemo

It is well known that Proposition 2.12 is false for locally compact groups in 
general (e.g., the motion groups). The following example 
shows that it can also fail for totally disconnected groups which are not 
compactly generated: 

\example{Example 2.13}Let $\Bbb Z^{*\Bbb N}=\{x\in \Bbb Z^{\Bbb N}\,;\, x_i 
\ne 0\text{ for finitely many }i\}$ and give $\Bbb Z^{*\Bbb N}$ the 
discrete topology. Give the multiplicative group $\{-1,1\}^{\Bbb N}$ the 
product topology. Let $\vp\:\{-1,1\}^{\Bbb N}\to\Aut(\Bbb Z^{*\Bbb N})$ be 
given by 
$\vp\bigl((\om_i)_{i=1}^\iy\bigr)\bigl((x_i)_{i=1}^\iy\bigr) 
=(\om_ix_i)_{i=1}^\iy$ and let $G$ be the semidirect product $G=\Bbb Z^{*\Bbb 
N}\ti\negmedspace_\vp\medspace\{-1,1\}^{\Bbb N}$. 

$\Bbb Z^{*\Bbb N}\ti\{e\}$ is trivially a closed cocompact normal SIN 
subgroup of $G$ but $G$ is not a SIN group because for every nonidentity 
element $g=(e,w)\in \{e\}\ti\{-1,1\}^{\Bbb N}$ there is $a\in \Bbb Z^{*\Bbb 
N}$ with $(a,e)(e,w)(a,e)^{-1}\notin \{e\}\ti\{-1,1\}^{\Bbb N}$. Indeed, if 
$w_j=-1$ and $a=(\del_{ji})_{i=1}^\iy$ then $(a,e)(e,w)(a,e)^{-1}=(v,w)$ 
where $v_j=2$. 
\endexample 

\proclaim{Lemma 2.14}Let $G$ be a locally compact compactly generated totally 
disconnected solvable group. Then there exists a closed normal 
cocompact subgroup $N$ such that $\ov{[G,G]}\sub N$ and $N/\ov{[G,G]}$ is 
topologically isomorphic to $\Bbb Z^d$ for some $d\ge 0$.
\endproclaim 
\demo{Proof}$G/\ov{[G,G]}$ is a compactly generated totally disconnected 
abelian group. Hence, it is the direct product $AB$ where $A\cong \Bbb Z^d$ 
and $B$ is a totally disconnected compact abelian group. Put $N=\pi^{-1}(A)$ 
where $\pi\:G\to G/\ov{[G,G]}$ is the canonical homomorphism. 
\qed\enddemo
\proclaim{Theorem 2.15}Conditions {\rm (i), (ii),} and {\rm (iii)} of 
Theorem 2.7 are equivalent when $G$ is a metrizable compactly 
generated totally disconnected locally compact metabelian group. 
\endproclaim 
\demo{Proof}Let $N$ be as in Lemma 2.14. To prove the nontrivial implication 
(ii)$\Ra$(iii) observe that as $\ov{[G,G]}$ is abelian, Theorem 2.7 applies 
to $N$. Thus if (ii) holds then $N$ is a SIN group. But then $G$ is a SIN 
group by Proposition 2.12. 
\qed\enddemo
\head{\bf 3. The Choquet-Deny theorem}\endhead

Let $\mu$ be a regular Borel probability measure on a locally compact group 
$G$. Recall that $G_\mu$ denotes the smallest closed subgroup containing the 
support of $\mu\/$. $\mu$ is called {\it adapted\/} if $G_\mu=G$. We will say 
that $\mu$ is a {\it Choquet-Deny measure\/} if every bounded continuous 
$\mu$-harmonic function is constant on the left cosets of $G_\mu$. 

We note that in the literature the Choquet-Deny theorem is often understood 
as the statement that every adapted $\mu\in M_1(G)$ is a Choquet-Deny measure 
(i.e., all bounded continuous $\mu$-harmonic functions are constant). We 
emphasize that in this paper the Choquet-Deny theorem is understood as the 
(formally) stronger statement that {\it every\/} $\mu\in M_1(G)$ is a 
Choquet-Deny measure. It is not known if the two versions of the Choquet-Deny 
theorem are equivalent. However, we know of examples of almost connected Lie 
groups with the property that every adapted spread out probability measure is 
Choquet-Deny but some non-adapted spread out measures are not. It can be 
shown (see Lemma 4.1) that the strong version of the theorem is true about 
$G$ if and only if the weak version holds for every closed subgroup of $G$. 

Throughout the sequel by the weak topology on the set $M_1(\Cal X)$ of 
probability measures on a locally compact space $\Cal X$ we mean the 
$\si(M_1(\Cal X),C_b(\Cal X))$-topology where $C_b(\Cal X)$ is the algebra of 
bounded continuous functions on $\Cal X$. 

\proclaim{Lemma 3.1} {\rm (a)} If $\mu$ is a Choquet-Deny measure on $G$ and 
$N\sub G$ is a closed normal subgroup, then the projection of $\mu$ onto 
$G/N$ is a Choquet-Deny measure on $G/N$. 

\noindent
{\rm (b)} If every neighbourhood of $e$ contains a compact normal 
subgroup $N$ such that the projection of $\mu$ onto $G/N$ is a 
Choquet-Deny measure, then $\mu$ is a Choquet-Deny measure. 
\endproclaim 
\demo{Proof}We omit a straightforward proof of (a). To prove (b) let 
us choose, for every neighbourhood $\Om$ of $e$, a compact normal 
subgroup $N_\Om\sub \Om$ such that the projection of $\mu$ onto 
$G/N_\Om$ is a Choquet-Deny measure. Denote by $\pi_\Om\:G\to G/N_\Om$ the 
canonical 
homomorphism and by $\om_\Om$ the normalized Haar measure of $N_\Om$. 
Directing the neighbourhoods of $e\/$ by reversed inclusion we obtain a 
net $(\om_\Om)$ in $M_1(G)$ which converges weakly to $\del_e$. 

Let $h$ be a bounded continuous $\mu$-harmonic 
function. We need to show that $h(xy)=h(x)$ for every $x\in G$ and $y\in 
G_\mu$. Now, when $\Om$ is a neighbourhood of $e$, the function $\om_\Om*h$ 
is a bounded continuous $\mu$-harmonic function constant on the cosets of 
$N_\Om$. Hence, $\om_\Om*h=h_\Om\circ\pi_\Om$ for a bounded continuous 
function $h_\Om$ on $G/N_\Om$. It is clear that $h_\Om$ is 
$\pi_\Om\mu$-harmonic where $\pi_\Om\mu$ denotes the 
projection of $\mu$ onto $G/N_\Om$. 
Moreover, $\pi_\Om(G_\mu)= (G/N_\Om)_{\pi_\Om\mu}$. 
Therefore for $x\in G$ and $y\in G_\mu$, 
$(\om_\Om*h)(xy)=h_\Om(\pi_\Om(x)\pi_\Om(y))=h_\Om(\pi_\Om(x)) 
=(\om_\Om*h)(x)$. 
Since $(\om_\Om*h)(\cdot)=\int_Gh(g^{-1}\cdot)\,\om_\Om(dg)$ 
and $\wlim_\Om\om_\Om=\del_e$, we conclude that $h(xy)=h(x)$. 
\qed\enddemo
\proclaim{Lemma 3.2}Let $(\mu_\al)$ be a net in $M_1(G)$. If for every 
neighbourhood $U$ of $e$ there exists $\vep\in M_1(G)$ such that $\vep(U)=1$ 
and the net $(\mu_\al *\vep)$ converges weakly, then the net $(\mu_\al)$ 
converges weakly.
\endproclaim 
\demo{Proof}There exists a compactly supported $\nu\in M_1(G)$ such 
that the net $(\mu_\al *\nu)$ is weakly convergent, and, hence, tight. This 
implies that the net $(\mu_\al)$ itself is tight. Then by Prohorov's theorem, 
every subnet of $(\mu_\al)$ has a weak cluster point. Therefore it suffices 
to show that the net $(\mu_\al)$ has a unique cluster point. But if $\mu'$ 
and $\mu''$ are cluster points of the net, then, due to our assumption, for 
every neighbourhood $U$ of $e$ there exists $\vep\in M_1(G)$ such that 
$\vep(U)=1$ and $\mu'*\vep=\mu''*\vep$. As in the proof of Lemma 3.1 we 
obtain a net $(\vep_i)$ in $M_1(G)$ which converges weakly to $\del_e$ and 
satisfies $\mu'*\vep_i=\mu''*\vep_i$ for every $i$. Hence, $\mu'=\mu''$. 
\qed\enddemo
\proclaim{Lemma 3.3}Let $G$ be a totally disconnected locally compact group, 
$\tau\in\Aut(G)$, and $F$ a finite subset of $C(\tau)$. If $\nu\in M_1(G)$ 
and $\nu(F)$=1 then the sequence $\nu*\tau\nu *\dots *\tau^{n-1}\nu$ 
converges weakly to a probability measure $\rho$ such that 
$\nu*\tau\rho=\rho$. 
\endproclaim 
\demo{Proof}It is clear that if $\rho=\wlim_{n\to\infty}\nu*\tau\nu *\dots 
*\tau^{n-1}\nu$ then $\nu *\tau\rho=\rho$. To see that the limit exists 
let $U$ be a compact open subgroup. Then there is $k\in\Bbb N$ 
such that for every $n\ge k$, $\tau^n(F)\sub U$. Let $\om_U$ denote the 
normalized Haar measure of $U$. Then for $n\ge k$, $\tau^n\nu*\om_U=\om_U$. 
Hence, $\nu*\tau\nu *\dots *\tau^{n-1}\nu*\om_U$ converges to $\nu*\tau\nu 
*\dots *\tau^{k-1}\nu*\om_U$. By Lemma 3.2, $\nu*\tau\nu *\dots 
*\tau^{n-1}\nu$ converges weakly. 
\qed\enddemo
\proclaim{Lemma 3.4}If $G$ is a locally compact group and $z\in G$ then 
$\ov{C(z)}\cap \<z\>=\{e\}$. 
\endproclaim 

\demo{Proof}Suppose that $z^k\in\ov{C(z)}$ for some $k>0$. Since 
$C(z)\sub C(z^k)$, we obtain $z^k\in\ov{C(z^k)}$. But when $U$ is a 
neighbourhood of $e$ in $\ov{C(z^k)}$, then 
$\ov{C(z^k)}=\ov{\bigcup_{n=1}^\iy 
z^{-kn}Uz^{kn}}$. This means that $\ov{C(z^k)}$ is either a {\it 
strange} group [19, Definition 1.1], or is compact. Since no locally compact 
group is strange [19, Theorem 1.8], $\ov{C(z^k)}$ is compact. As 
$z^k\in\ov{C(z^k)}$, it follows that $C(z)=\{e\}$. 
\qed\enddemo
Suppose that the locally compact group $G$ acts on a locally compact space 
$\Cal X$ so that the mapping $G\ti\Cal X\ni (g,x)\to gx\in\Cal X$ is 
continuous. Given $\rho\in M_1(\Cal X)$ and $g\in G$ we write $g\rho$ for 
the measure $(g\rho)(\cdot)=\rho(g^{-1}\cdot)$. Given $\mu\in M_1(G)$ we 
denote by $\mu *\rho$ the measure $(\mu*\rho)(\cdot) =\int_G 
(g\rho)(\cdot)\,\mu(dg)$. Now, if $\rho=\mu*\rho$ then for every bounded 
continuous function $f\:\Cal X\to \Bbb C$, the function $h(g)=\int_{\Cal 
X}f(gx)\,\rho(dx)=\int_{\Cal X}f(x)\,\,(g\rho)(dx)$ is a bounded continuous 
$\mu$-harmonic function. Therefore in order to show that the Choquet-Deny 
theorem fails for $\mu$ it suffices to find $g\in G_\mu$ such that 
$g\rho\ne\rho$. This observation is being used in the proof of the next 
lemma. 
\proclaim{Lemma 3.5} Let $G$ be a totally disconnected locally compact group 
and $z$ an element of $G$ with $C(z)\ne\{e\}$. Let $g\in C(z)-\{e\}$, and 
$\nu=p\del_g+(1-p)\del_e$ where $p\in (0,1)-\{\frac12\}$. Then the 
Choquet-Deny theorem is false for the measure $\mu=\nu*\del_z$. 
\endproclaim 
\demo{Proof}Let $\tau$ denote the the inner automorphism $z\cdot z^{-1}$. By 
Lemma 3.3 the limit $\rho=\wlim_{n\to\iy}\nu *\tau\nu *\dots *\tau^{n-1}\nu$ 
exists and satisfies $\nu*\tau\rho=\rho$. Moreover, $\rho(\ov{C(z)})=1$. 

Note that $\<z\>$ is necessarily infinite and discrete, so it is a closed 
subgroup of $G$ (isomorphic to $\Bbb Z$). Let $\pi\:G\to G/\<z\>$ denote the 
canonical mapping and let $\hat\rho =\pi\rho$. Then 
$\mu*\hat\rho=\pi(\mu*\rho) =\pi(\nu*\del_z*\rho)=\pi(\nu*\tau\rho 
*\del_z)=\pi\rho=\hat\rho$. Since $g\in G_\mu$, it suffices to show that 
$g\hat\rho\ne\hat\rho$. 

Now, there exists a compact subgroup $U$ of $\ov{C(z)}$ such that 
$g\nin U$ but $\tau^j(g)\in U$ for every $j\ge 1$. Let $\om_U$ be the 
normalized Haar measure of $U$. Then\lb $\nu *\tau\nu *\dots 
*\tau^{n-1}\nu*\om_U=\nu*\om_U=p(g\om_U)+(1-p)\om_U$. Thus 
$\rho*\om_U=p(g\om_U)+(1-p)\om_U$ and 
$g(\rho*\om_U)=p(g^2\om_U)+(1-p)(g\om_U)$. 
Since $\rho(\ov{C(z)})=1$ and by Lemma 3.4 $\ov{C(z)}\cap\<z\>=\{e\}$, we 
obtain $\hat\rho(\pi(U))=\rho(U\<z\>)=\rho(U)=(\rho*\om_U)(U)=1-p$ and 
$(g\hat\rho)(\pi(U))=(g\rho)(U\<z\>)=(g\rho)(U)= 
(g(\rho*\om_U))(U)=p(g^2\om_U)(U)\ne 1-p$. 
\qed\enddemo

\proclaim{Theorem 3.6} Let $G$ be a totally disconnected generalized 
$\ov{FC}$-group or a metrizable locally compact compactly generated totally 
disconnected metabelian group. Then the following conditions are 
equivalent\rom:
\roster
\item"{(a)}" The Choquet-Deny theorem holds for $G$. 
\item"{(b)}" The Choquet-Deny theorem holds for every $\mu\in M_1(G)$ with 
             $\supp\mu$ of cardinality 2.
\item"{(c)}" $G$ is distal and has polynomial growth. 
\endroster
\endproclaim 
\demo{Proof}(b)$\Ra$(c): We first prove that $G$ is distal. When $G$ is 
metrizable, this is clear by Lemma 3.5, Theorems 2.9 and 2.15, and 
Proposition 2.1. Suppose that $G$ is a not necessarily metrizable generalized 
$\ov{FC}$-group. Note that it suffices to show that every neighbourhood of 
$e\/$ contains a compact normal subgroup $N$ such that $G/N$ is distal. But 
as 
$G$ is compactly generated, given a neighbourhood $U$ of $e$ there exists a 
compact normal subgroup $N\sub U$ such that $G/N$ is metrizable [8, Theorem 
8.7]. Every probability measure on $G/N$ with support of cardinality 2 is the 
canonical image of a similar measure on $G$. Hence, by Lemma 3.1(a), 
Condition (b) must hold on $G/N$ and as $G/N$ is a generalized 
$\ov{FC}$-group, it is distal. 

We now prove that $G$ is of polynomial growth. Suppose that $G$ is not of 
polynomial growth. By Proposition 2.8 and Theorem 2.15, $G$ has a compact 
open normal subgroup $N$ such that the quotient $G/N$ contains a finitely 
generated solvable subgroup $S$ of finite index (polycyclic when $G$ is an 
$\ov{FC}$-group and metabelian when $G$ is metabelian). By [10, Th\'eor\`eme 
I.4] $S$ is not of polynomial growth. Hence, by [18, Theorem 3.13 and its 
proof], $S$ supports a probability measure with a 
2-element support for which the Choquet-Deny theorem fails. This implies that 
the Choquet-Deny theorem fails for a similar probability measure on $G$. 

(c)$\Ra$(a): When $N$ is a compact open normal subgroup, $G/N$ is a finitely 
generated group of polynomial growth, hence, the Choquet-Deny theorem holds 
for $G/N$. Since by Theorems 2.9 and 2.15, $G$ has arbitrarily small compact 
open normal subgroups, Lemma 3.1(b) yields the desired conclusion. 
\qed\enddemo
\head{\bf 4. On boundaries of random walks} \endhead 

It is well known that the bounded $\mu$-harmonic functions can be 
represented, by means of a ``Poisson formula", as bounded Borel functions on 
a certain 
``boundary space". Let us consider the bounded $\mu$-harmonic functions (on a 
general locally compact group $G$) as elements of $L^\iy(G)$, and let 
$\Cal H_\mu$ denote the resulting subspace of $L^\iy(G)$. $\Cal H_\mu$ is 
invariant under the usual left action of $G$ on $L^\iy(G)$ and for every 
absolutely continuous $\nu\in M_1(G)$ and every $h\in\Cal H_\mu$, $\nu*h$ is 
a bounded (left uniformly) continuous $\mu$-harmonic function. When $G$ is 
locally compact second countable (lcsc), there exists a standard Borel 
$G$-space $\Cal X$ with a $\si$-finite quasiinvariant measure $\al$ and 
an equivariant isometry $\Phi$ of $L^\iy(\Cal X,\al)$ onto $\Cal H_\mu$ 
[15, \S 3]. $\Phi$ is given by the Poisson formula 
$$
(\Phi f)(g)=\int_{\Cal X}f(gx)\,\rho (dx)                  \tag{4.1}
$$
where $\rho$ is a probability measure on $\Cal X$ satisfying $\mu*\rho=\rho$.
The $G$-space $\Cal X\/$, called the {\it $\mu$-boundary\/}, or {\it Poisson 
boundary\/}, is not unique. However, for any two $\mu$-boundaries $(\Cal 
X',\al')$ and $(\Cal X'',\al'')$, there exists an equivariant isomorphism 
between $L^\iy(\Cal X',\al')$ and $L^\iy(\Cal X'',\al'')$ (which implies that 
$(\Cal X',\al')$ and $(\Cal X'',\al'')$ are isomorphic up to sets of zero 
measure). The $\mu$-boundary can be always realized as a topological, compact 
metric $G$-space [15, \S 3]. 

When the Choquet-Deny theorem holds for $\mu$, the natural realization of the 
$\mu$-boundary is the homogeneous space $G/G_\mu$ where the ``Poisson kernel" 
$\rho$ (cf. Eq. 4.1) is the point measure $\del_{G_\mu}$. When $G$ is a 
discrete (countable) group then the $\mu$-boundary is a homogeneous space if 
and only if the Choquet-Deny theorem is true for $\mu$ [18, Lemma 1.1 and 
the remark preceding Proposition 2.6]. The situation is different for 
continuous groups. When $G$ is 
an almost connected lcsc group then for every spread out probability measure 
on $G$ the $\mu$-boundary is a homogeneous space [14, Corollary 4.7].
It is well known that the $\mu$-boundary of every spread out 
measure on a connected semisimple Lie group with finite centre is a compact 
homogeneous space [6,2]. However, if the $\mu$-boundary of 
a spread out measure on an {\it amenable\/} lcsc group is a {\it compact\/} 
homogeneous space, then the Choquet-Deny theorem holds for 
$\mu$ (and the $\mu$-boundary is finite), see [18, Proposition 2.6] and 
[16, Lemma 2.3], or [2, Propositions IV.8 and IV.7]. 

Theorem 4.2 which we prove below applies, in particular, to every totally 
disconnected compactly generated lcsc group of polynomial growth. The 
result is that for such groups the $\mu$-boundary can be always realized as a 
homogeneous space, and, as a compact homogeneous space when $\mu$ is adapted; 
when $\mu$ is adapted and spread out the $\mu$-boundary is a singleton. 

\proclaim{Lemma 4.1}A probability measure $\mu$ on a locally compact group 
$G$ is a Choquet-Deny measure if and only if the restriction of $\mu$ 
to $G_\mu$ is a Choquet-Deny measure (on $G_\mu$). 
\endproclaim 
\demo{Proof}Let $\mu'$ denote the restriction of $\mu$ to $G_\mu$. The 
restriction of a $\mu$-harmonic function to $G_\mu$ is $\mu'$-harmonic; 
moreover, if $h$ is $\mu$-harmonic then for every $g\in G$ the left translate 
$(gh)(\cdot)=h(g^{-1}\cdot)$ is also $\mu$-harmonic. Hence, if $\mu'$ is 
Choquet-Deny then so is $\mu$. The converse is equally obvious 
when $G_\mu$ is open, because then every bounded continuous $\mu'$-harmonic 
function trivially extends to a bounded continuous $\mu$-harmonic function. 
However, in general, a technical argument is called for. 

Let us first consider the case that $G$ is second countable. Let $h'$ be a 
bounded continuous $\mu'$-harmonic function. As $G$ is second countable, 
the canonical projection $\pi\:G\to G/G_\mu$ admits a Borel cross-section 
$\k$. Since for every $g\in G$, $\k(\pi(g))^{-1}g\in G_\mu$, we can 
define a function $h\:G\to \Bbb C$ by $h(g)=h'\bigl(\k(\pi(g))^{-1}g\bigr)$. 
$h\phantom{\ }$is a bounded (in general, discontinuous) $\mu$-harmonic 
function. 

Let $(\vep_n)$ be a sequence of absolutely continuous probability measures on 
$G$ converging weakly to $\del_e$. Then the sequence $(\vep_n*h)$ converges 
in the weak* topology of $L^\iy(G)$ to $h$. Since $\vep_n*h$ is a bounded 
continuous $\mu$-harmonic function and $\mu$ is a Choquet-Deny measure, it 
follows that there exists a bounded Borel function $\hat h\:G/G_\mu\to\Bbb C$ 
such that $h=\hat h\circ\pi$ $\lam$-a.e., where $\lam$ is the Haar measure of 
$G$. 

Now, the mapping $\vp\:(G/G_\mu)\ti G_\mu\to G$ given by $\vp(x,g)=\k(x)g$ is 
a Borel isomorphism. Moreover, if $\nu$ is a $\si$-finite quasiinvariant 
measure on $G/G_\mu$ and $\lam'$ the Haar measure of $G_\mu$, then the 
measure $\vp(\nu\ti\lam')=(\nu\ti\lam')\circ\vp^{-1}$ is equivalent to the 
Haar measure $\lam$ of $G$. Consequently,
$$
\gather
0=\inl_{(G/G_\mu)\ti G_\mu}|h'\circ\vp-\hat h\circ\pi\circ\vp 
|\,\,d(\nu\ti\lam')
\\
=\inl_{G/G_\mu}\biggl[\,\inl_{G_\mu}|(h\circ\vp)(x,g)-(\hat 
h\circ\pi\circ\vp)(x,g)|\,\,\lam'(dg)\biggr]\,\nu(dx)
\\
=\inl_{G/G_\mu}\biggl[\,\inl_{G_\mu}|h'(g)-\hat h(x)|\,\,\lam'(dg)\biggr] 
\,\nu(dx).  
\endgather
$$
Thus for $\nu$-a.e. $x\in G/G_\mu$, $\int_{G_\mu}|h'(g)-\hat h(x)|\,\lam'(dg)
=0$. Hence, as $h'$ is continuous, it is constant. 

Consider now the general case that $G$ is not necessarily second countable. 
Observe that due to the regularity of $\mu$ and local compactness of 
$G$, $G_\mu$ is $\si$-compact and, hence, there is also an open $\si$-compact 
subgroup $G_1$ with $\mu(G_1)=1$. Since $G_1$ is open it is clear that the 
restriction of $\mu$ to $G_1$ is a Choquet-Deny measure. Hence, we may assume 
that $G$ itself is $\si$-compact. By Lemma 3.1(b) it suffices to show that 
every neighbourhood $U$ of $e$ 
in $G_\mu$ contains a compact normal subgroup $N$ such that the projection of 
$\mu'$ onto $G_\mu/N$ is Choquet-Deny. But by [8, Theorem 8.7] there exists a 
compact normal subgroup $K$ of $G$ such that $K\cap G_\mu\sub U$ and $G/K$ is 
second countable. Let $\pi_K\:G\to G/K$ denote the canonical homomorphism. 
Since $(G/K)_{\pi_K\mu}=\pi_K(G_\mu)$, combining Lemma 3.1(a) with what we 
just proved~for second countable groups, we conclude that the restriction of 
$\pi_K\mu$ to $\pi_K(G_\mu)$ is a Choquet-Deny measure. As $\pi_K(G_\mu)$ is 
canonically isomorphic to $G_\mu/(K\cap G_\mu)$, it follows that the 
projection of $\mu'$ onto $G_\mu/(K\cap G_\mu)$ is a Choquet-Deny measure.\qed
\enddemo

\proclaim{Theorem 4.2}Let $\mu\in M_1(G)$ where $G$ is a lcsc group. If $G$ 
contains a compact normal subgroup $K$ such that the projection of $\mu$ onto 
$G/K$ is a Choquet-Deny measure, then the $\mu$-boundary can be realized as a 
homogeneous space; when $\mu$ is adapted, the $\mu$-boundary can be realized 
as a compact homogeneous space on which $K$ acts transitively. 
\endproclaim 
\demo{Proof}Denote by $\pi\:G\to G/K$ the canonical homomorphism.

Suppose that $\mu$ is adapted and let $(\Cal X,\al)$ be the $\mu$-boundary 
realized as a standard Borel $G$-space.  Let $f\in L^\iy(\Cal X,\al)$ be 
invariant under the action of $K$. Then the corresponding $\mu$-harmonic 
function $h=\Phi f\in \Cal H_\mu$ (cf. Eq. 4.1) is also invariant under the 
(left) action of $K$. Hence, $h=\hat h\circ\pi$ where $\hat h\in\Cal 
H_{\pi\mu}$. Since $\pi\mu$ is adapted and the Choquet-Deny theorem holds on 
$G/K$, it follows that $h$ is constant. Thus so is $f$. 
As $\Cal X$ is a standard Borel 
$G$-space this implies that $K$ acts ergodically on $\Cal X$, and, hence, 
$\al$ is carried on an orbit of $K$ [30, Corollary 2.1.21 and Proposition 
2.1.10]. Consequently, the $\mu$-boundary can be realized as a compact 
homogeneous space of $G$ on which $K$ acts transitively. 

When $\mu$ is not necessarily adapted, let $\mu'$ denote the restriction of 
$\mu$ to $G_\mu$ and let $(\Cal X',\al')$ be a realization of the 
$\mu'$-boundary as a standard Borel $G_\mu$-space. By Lemma 4.1 the 
restriction of $\pi\mu$ to $\pi(G_\mu)=(G/K)_{\pi\mu}$ is a Choquet-Deny 
measure. Since $G_\mu/(G_\mu\cap K)\cong \pi(G_\mu)$, it follows that the 
projection of $\mu'$ onto $G_\mu/(G_\mu\cap K)$ is a Choquet-Deny measure. As 
$\mu'$ is adapted, we may assume that $\Cal X'$ is a homogeneous space of 
$G_\mu$ (on which $G_\mu\cap K$ acts transitively). 
Now, by [17, Proposition 3.5 and Remark 3.9] the $\mu$-boundary can be 
realized as the skew product $\Cal X=G/G_\mu\ti\negmedspace_\g\medspace\Cal 
X'$ (the $G$-space induced from the $G_\mu$-space $\Cal X'$ [30, p.\,75]), 
where $\g\:G\ti G/G_\mu\to G_\mu$ is the cocycle associated with a Borel 
cross section of the canonical projection of $G$ on $G/G_\mu$. It follows 
that $G$ acts transitively on $\Cal X$. This means that the $\mu$-boundary 
can be realized as a homogeneous space of $G$.
\qed\enddemo
\proclaim{Corollary 4.3}Let $G$ be a totally disconnected compactly generated 
lcsc group of polynomial growth. 
Then for every $\mu\in M_1(G)$ the 
$\mu$-boundary can be realized as a homogeneous space of $G$; when $\mu$ is 
adapted, the $\mu$-boundary can be realized as a compact homogeneous space. 
\endproclaim 
The next corollary can be regarded as a generalization of the implication 
(c)$\Ra$(a) of Theorem 3.6. Contrary to the proof of Theorem 3.6, the proof 
of Corollary 4.4 does not rely on equicontinuity of $\Inn(G)$. 

\proclaim{Corollary 4.4}Let $G$ be a locally compact group containing a 
compact normal subgroup $K$ such that the Choquet-Deny theorem holds for 
$G/K$ and $\,\Inn(G)$ acts distally on $K$. Then the Choquet-Deny theorem 
holds for $G$. 
\endproclaim
\demo{Proof}It is not difficult to see that if a locally compact group $G$ 
contains a compact normal subgroup $K$ such that the Choquet-Deny theorem 
holds for $G/K$ and $\Inn(G)$ acts distally on $K$, then the same is true for 
every closed subgroup and every quotient of $G$. Let $\mu\in M_1(G)$. To show 
that $\mu$ is a Choquet-Deny measure it suffices to show that the 
restriction, $\mu'\/$, of $\mu$ to $G_\mu$ is Choquet-Deny. By Lemma 3.1(b), 
to show the 
latter it is enough to show that every neighbourhood of $e$ in $G_\mu$ 
contains a compact normal subgroup $N$ such that 
the projection of $\mu'$ onto $G_\mu/N$ is Choquet-Deny. But as $G_\mu$ 
is $\si$-compact, every neighbourhood of $e$ in $G_\mu$ contains a compact 
normal subgroup with second countable quotient. 
Hence, it is enough to prove that if a 
lcsc group $G$ contains a compact normal subgroup $K$ such that the 
Choquet-Deny theorem holds for $G/K$ and $\Inn(G)$ acts distally on $K$, then 
every adapted probability measure on $G$ is a Choquet-Deny measure. 

For such $G$ and $\mu$, by Theorem 4.2, the $\mu$-boundary has the 
form $G/H$ where $K$ acts transitively on $G/H$, i.e., $G=KH$. Let $\rho$ 
denote the Poisson kernel. Note that due to the identity $\rho=\mu*\rho$ and 
adaptedness of $\mu$, it suffices to show that $\rho$ is a point measure 
(this will imply that $G/H$ is a singleton). 

Now, by [11, Proposition 1.8] there exists a sequence $(h_n)$ in $G$ such 
that the sequence $(h_n\rho)$ converges weakly to a point measure 
$\del_{x_0}$. Since $G=KH$ and $K$ is compact, we may assume that $h_n\in H$ 
for all $n$. Next, by [19, Lemma 2.8] we may assume that there is a Borel set 
$B\sub G/H$ such that $\rho(B)=1$ and $\lim_{n\to\infty}h_nx =x_0$ for every 
$x\in B$. It is enough to show that $B$ is a singleton. 

Consider the compact homogeneous space $K/(K\cap H)$. The formulas 
$h^\centerdot k=hkh^{-1}$ and $h_\centerdot k(K\cap H)=hkh^{-1}(K\cap H)$, 
$k\in K$, define actions of $H$ on $K$ and $K/(K\cap H)$, respectively. 
Clearly, the $\vphantom h_\centerdot\,$-action is a factor of the 
$\vphantom h^\centerdot\,$-action. As $\vphantom h^\centerdot$ is distal, so 
is $\vphantom h_\centerdot$ [5 , Corollary 6.10, p.\,52]. 

Let $x_1,x_2\in B$. Write $x_1=k_1H$ and $x_2=k_2H$ with $k_1,k_2\in K$. 
Then $\lim_{n\to\infty}h_nx_j =\lim_{n\to\infty}h_nk_jh_n^{-1}H =x_0$ for 
$j=1,2$. Since the compact homogeneous spaces $G/H$ and $K/(K\cap H)$ are 
isomorphic as $K$-spaces, we then obtain $\lim_{n\to\infty}h_n\vphantom h_
\centerdot k_1(K\cap H)=\lim_{n\to\infty}h_n\vphantom h_\cdot k_2(K\cap H)$. 
Since $\vphantom h_\centerdot$ is distal, $k_1(K\cap H)=k_2(K\cap H)$ and, 
hence, $x_1=x_2$. Therefore $B$ is singleton. 
\qed\enddemo
\example{Example 4.5}Let $\tau$ be the automorphism of the torus $\Bbb T^3$, 
defined by $\tau(x,y,z)=(x,xy,xyz)$. Then $\tau$ is distal but not 
equicontinuous. By Corollary 4.4 the Choquet-Deny theorem is true for the 
3-step nilpotent group $\Bbb T^3\ti\negmedspace_\tau\medspace\Bbb Z$. 
\endexample 
\example{Example 4.6}Let $\tau$ be the shift $\tau((x_i)_{i\in\Bbb Z})=
(x_{i+1})_{i\in\Bbb Z}$ on the compact abelian group $\Bbb Z_2^{\Bbb Z}$, 
and let $G=\Bbb Z_2^{\Bbb Z}\ti\negmedspace_\tau\medspace\Bbb Z$. 
Since $C(\tau)=\{x\in \Bbb Z_2^{\Bbb Z}\,;\,\text{there exists }
k\nmb\in\nmb\Bbb Z\text{ with } x_i=0\text{ for every }i\ge k\}$, 
$\Inn(G)$ does not act distally on $\Bbb Z_2^{\Bbb Z}\ti\{0\}$. By Theorem 
3.6 the Choquet-Deny theorem is not true for $G$. 

Let $\mu\in M_1(G)$ be adapted. According to Theorem 4.2 the $\mu$-boundary 
has the form $G/H$ where $G=(\Bbb Z_2^{\Bbb Z}\ti\{0\})H$. It is not 
difficult to see that a closed subgroup $H\sub G$ satisfies $G=(\Bbb 
Z_2^{\Bbb Z}\ti\{0\})H\/$ if and only if there is a closed $\tau$-invariant 
subgroup $T\sub\Bbb Z_2^{\Bbb Z}$ and $g\in G$ such that $gHg^{-1}=T\ti\Bbb 
Z$. We may therefore assume that $H=T\ti\Bbb Z$ where $T$ is a closed 
$\tau$-invariant subgroup. Now, the formula $(x,y)(zT)=x\tau^y(z)T$, $x,z\in 
\Bbb Z_2^{\Bbb Z}$, $y\in\Bbb Z$, defines an action of $G$ on $\Bbb Z_2^{\Bbb 
Z}/T$ under which $\Bbb Z_2^{\Bbb Z}/T$ becomes a homogeneous space of $G$, 
isomorphic to $G/H$. Thus for an adapted $\mu\in M_1(G)$, the $\mu$-boundary 
can be realized as one of the $G$-spaces 
$\Bbb Z_2^{\Bbb Z}/T$, where $T$ is a closed $\tau$-invariant subgroup of 
$\Bbb Z_2^{\Bbb Z}$. 

Let for $k=1,2,\dots$, $\,S_k=\{x\in\Bbb Z_2^{\Bbb Z}\,\;\, \tau^k(x)=x\}$. 
Then $S_k$ is a closed $\tau$-invariant subgroup of $\Bbb Z_2^{\Bbb Z}$ and 
it can be shown that $T$ is a closed $\tau$-invariant subgroup of $\Bbb 
Z_2^{\Bbb Z}$ if and only if $T=\Bbb Z_2^{\Bbb Z}$ or $T$ is a 
$\tau$-invariant subgroup of $S_k$ for some $k$. In particular, proper 
$\tau$-invariant subgroups are finite. Let $\Cal T$ denote the class of 
closed $\tau$-invariant subgroups of $\Bbb Z_2^{\Bbb Z}$. The $G$-spaces 
$\Bbb Z_2^{\Bbb Z}/T$, $T\in\Cal T$ are mutually nonisomorphic and each of 
them is an equivariant image of $\Bbb Z_2^{\Bbb Z}=\Bbb Z_2^{\Bbb Z}/S_1$. 

One can construct a family $\mu_T$, $T\in\Cal T$, of discrete probability 
measures on $G$ such that $\Bbb Z_2^{\Bbb Z}/T$ is the $\mu_T$-boundary for 
every $T\in\Cal T$. We refrain from going into the details here as this 
would require a longer digression into the theory of the $\mu$-boundaries. A 
more difficult question concerns determining, when $\mu\in M_1(G)$ 
is given, which of the spaces $\Bbb Z_2^{\Bbb Z}/T$ is the $\mu$-boundary. 
In particular, one would like to know for which $\mu\in M_1(G)$ the 
$\mu$-boundary is a singleton. In addition to the case of spread out 
measures, this is so for every adapted probability measure which induces a 
recurrent random walk on $\Bbb Z\cong G/(\Bbb Z^{\Bbb Z}_2\ti\{0\})$. We do 
not know of any relevant conditions that are both sufficient and necessary. 
\endexample 

\phantom{oooooooooooooooooo}

\flushpar
{\smc Acknowledgement} 

\flushpar
The second author thanks the School of Mathematics and Statistics of Carleton 
University for its hospitality.

\phantom{oooooooooooooooooo}

\enddocument